\documentclass{article}%
\usepackage{amsmath}
\usepackage{amsfonts}
\usepackage{amssymb}
\usepackage{graphicx}%
\begin{document}

\title{A Method for Preserving Geometric Meaning in the Calculus of Variations}
\author{Pavel Grinfeld\\Drexel University}
\maketitle

\begin{abstract}
We discuss the advantages of the Calculus of Moving Surfaces over the
Euler-Lagrange equation for optimization problems originating in Geometry. An
extension of Tensor Calculus, it provides tools for analyzing geometric
quantities directly rather than their coordinate representations. This allows
us to avoid the many difficulties associated with the use of coordinates, from
the untenable complexity of analytical expressions to the virtual
impossibility of recovering the geometric interpretation of the final result.
As an illustration, we analyze the brachistochrone and give its geometric
characterization in terms of its curvature.

\end{abstract}

\section{Introduction}

Every geometric problem in the Calculus of Variations is a problem in the
Calculus of Moving Surfaces and therefore admits a more elegant and
geometrically revealing solution. We will make this evident by considering the brachistochrone.

In $1744$ Leonhard Euler published his celebrated \textit{Method for Finding
Curves Enjoying the Property of Maximum or Minimum \cite{Euler1744}}, in which
he laid the foundations of the Calculus of Variations. From its very first
pages, we see a strong drive toward the algebraization of Geometry by the
method of coordinates. In Paragraph 4 of Chapter I we read \textit{...the
precise study of curves requires that they be referred to some fixed
axis...}The trend thus set by Euler continues to this day.

Joseph Louis Lagrange was inspired by this approach and took it further in his
$1761$ \textit{Essay on a New Method for Finding the Maxima And Minima of
Indefinite Integral Expressions \cite{LagrangeCalculusOfVariations}}.\ In this
essay, he reduced the Calculus of Variations to a purely analytical method
that culminated in the Euler-Lagrange equations. Lagrange's approach greatly
broadened the range of problems that can be solved and, in particular,
extended the analysis from curves to surfaces. In one of his letters to
Lagrange, Euler wrote: \textit{Your analytical solution of the isoperimetric
problem, from what I\ can see, leaves nothing to be desired in this field of
inquiry, and I am delighted beyond measure that it has been your good fortune
to carry this subject -- which, from its inception, I had cultivated almost
entirely on my own -- to such a high degree of perfection.}

We believe, however, that Lagrange's method can be taken to an even higher
level of perfection. The remarkable effectiveness of the Euler-Lagrange
equations came at the expense of geometric meaning. Furthermore, the same
observation applies to the method of coordinates in general. Untethering an
analysis from the underlying geometry invites a cascade of woes. First is the
disorienting flatness of the analytical expressions that deprives us of
logical thrust -- the kind of thrust that is found in all great mathematical
analyses. Second is the growing analytical complexity that often becomes
untenable and forces us to retreat in the face of computational difficulties.
Lastly, if we press forward and, improbably, succeed in reaching a final
answer, we are unlikely to know its simple geometric meaning, thus preventing
it from becoming a discovery.

The very works of Euler and Lagrange bear witness to this fact. In Chapter
$5$, Paragraph $47$ of his \textit{Method}, Euler, with his remarkable ability
to find the simplest example that illustrates the general problem, solved for
the minimal surface of revolution. He showed that the desired profile
$y\left(  x\right)  $ satisfies the equation%
\begin{equation}
y_{xx}y-y_{x}^{2}-1=0. \label{Euler}%
\end{equation}
Euler promptly recognized that the solution is the famous \textit{catenary},
i.e. the shape of a hanging chain.

Tellingly, Euler did not find -- nor, in truth, did he care to find -- the
geometric interpretation of this equation, which is that the \textit{mean
curvature} $B_{\alpha}^{\alpha}$, i.e. the sum of the principal curvatures,
vanishes:%
\begin{equation}
B_{\alpha}^{\alpha}=0. \label{BAa = 0}%
\end{equation}
It should be noted, however, that the concept of principal curvatures would be
introduced by Euler himself in $1760$ \cite{Euler1767Echerches}, and was
therefore not available to him in $1744$.

A similar absence of geometric interpretation is found in Lagrange's work.
Lagrange showed that a minimal surface represented by a function $z\left(
x,y\right)  $ is characterized by the equation%
\begin{equation}
\left(  1+z_{y}^{2}\right)  z_{xx}+2z_{x}z_{y}z_{xy}+\left(  1+z_{x}%
^{2}\right)  z_{yy}=0, \label{Lagrange}%
\end{equation}
and, as in Euler's work, no geometric interpretation was given.

The connection between minimal surfaces and mean curvature was found in $1776$
by the French mathematician Jean Baptiste Meusnier \cite{Meusnier1785Memoire}.

Future generations of mathematicians have sought to restore the connection to
the underlying geometry in a number of ways. One such avenue led to the
development of Tensor Calculus by Gregorio Ricci and Tullio Levi-Civita at the
turn of the twentieth century. However, Tensor Calculus is not sufficient for
the Calculus of Variations as the latter considers families of surfaces, the
analysis of which comes with its peculiar artifacts in the coordinate
approach. Thus, the Calculus of Variations requires the \textit{Calculus of
Moving Surfaces} which will prove to be of a "higher degree of perfection"
than the conventional approach.

This is especially true for variational problems originating in Geometry, for
there is a mismatch between how those problems are formulated and the form
required by the conventional approach to the Calculus of Variations. When a
geometric quantity is to be minimized or maximized, possibly subject to one or
more constraints, that quantity and the constraints are usually expressed in
the form of a \textit{geometric} integral, such as the integral%
\begin{equation}
\int_{S}dS
\end{equation}
for the total area of a surface. (The term \textit{geometric} in this context
means that the integral is defined in purely geometric terms. Informally, the
integral%
\begin{equation}
\int_{S}FdS
\end{equation}
captures the \textit{total amount} of $F$ over the surface. A rigorous
definition is not needed here.) Meanwhile, the Euler-Lagrange equation
requires that same quantity to be expressed as an \textit{arithmetic} integral%
\begin{equation}
\idotsint F\left(  x^{i},z^{\alpha},\frac{\partial z^{\alpha}}{\partial x^{i}%
}\right)  dx^{1}\cdots dx^{n},
\end{equation}
such as%
\begin{equation}
\int\int\sqrt{1+z_{x}^{2}+z_{y}^{2}}dxdy
\end{equation}
for the area of a surface. Therefore, the Euler-Lagrange equation effectively
requires that the problem be referred to a specific coordinate system which,
as we have described, relinquishes much of the connection to the original
geometry along with almost any chance of recovering it in the end.

By contrast, the Calculus of Moving Surfaces can be applied directly to the
geometric formulation of the problem. Being an extension of Tensor Calculus,
the Calculus of Moving Surfaces does, at its core, rely on the use of
coordinates, but combines the analytical expressions in a way that allows for
the preservation of the geometric meaning.

In what follows, we will see this promise in action.

\section{The moving surfaces formulation of variational problems}

Suppose that for a given Euclidean manifold $S$ -- be it a curve, a surface,
or a volume -- we have an optimization problem in which the objective quantity
$E$ is expressed in terms of line, surface, or volume integrals, while the
independent variation is the hypersurface $S$ itself, i.e. its shape and
location. We will limit our attention to hypersurfaces, i.e. curves embedded
in a plane or a surface, and surfaces embedded in space.

As a model problem, consider the classical isoperimetric problem, which is to
minimize the total surface area%
\begin{equation}
E=\int_{S}dS
\end{equation}
of a closed surface $S$ subject to the enclosed volume constraint%
\begin{equation}
\int_{\Omega}d\Omega=V_{0},
\end{equation}
where $\Omega$ is the enclosed domain.

Let us proceed by imitating the key analytical device in the derivation of the
Euler-Lagrange equations. That is, we will consider a family of hypersurfaces
$S\left(  t\right)  $ such that the optimal surface occurs at $t=0$. We can
intuitively think of $t$ as time and of $S\left(  t\right)  $ as a
\textit{moving surface}. Naturally, the quantity $E$ can be evaluated at any
$t$ and therefore can be treated as a function of time. For example,%
\begin{equation}
E\left(  t\right)  =\int_{S\left(  t\right)  }dS.
\end{equation}
Likewise, the constraint also becomes dependent on time, i.e.%
\begin{equation}
\int_{V\left(  t\right)  }dV=V_{0}.
\end{equation}
Then for \textit{any} family $S\left(  t\right)  $, the function $E\left(
t\right)  $ attains its extremal value at $t=0$ and therefore%
\begin{equation}
E^{\prime}\left(  0\right)  =0.
\end{equation}
This is the very condition from which, by taking advantage of the
arbitrariness of $S\left(  t\right)  $, we can derive the equilibrium
equations. What is required to turn this blueprint into a practical procedure
is the machinery of the Calculus of Moving Surfaces.

\section{The essential elements of the\ Calculus of Moving Surfaces}

For the sake of brevity, we will merely name, rather than define or discuss,
the essential elements of the Calculus of Moving Surfaces required for our
presentation. For a thorough introduction to the subject, we refer the reader
to the paper \textit{A Better Calculus of Moving Surfaces} \cite{BetterCMS},
which gives a brief historical overview of the subject and defines the key
differential operator $\dot{\nabla}$ in its present form. Meanwhile, a
detailed exposition of the subject can be found in \cite{GrinfeldTC}.

The motion of the hypersurface $S\left(  t\right)  $ is characterized by the
invariant \textit{normal velocity} $C$. The scalar field $C$ can be
intuitively understood as the rate of deformation of the surface in the normal
direction. In particular, the sign of $C$ depends on the chosen direction of
the \textit{unit normal} $\mathbf{N}$. Furthermore, the combination
$C\mathbf{N}$, which is independent of the choice of normal, can be thought of
as the \textit{vector} rate of deformation in the normal direction.

The symbol $\dot{\nabla}$ represents the \textit{invariant time derivative
with respect to the moving surface}. When applied to an invariant, such as
$C$, $\mathbf{N}$, or $B_{\alpha}^{\alpha}$ (or any non-invariant tensor of
order zero), it is interpreted as the rate of change of that quantity along
the normal direction. The analytical definition of $\dot{\nabla}$ is quite
elaborate and can be found in \cite{BetterCMS} and \cite{GrinfeldTC}. The one
property of $\dot{\nabla}$ that is relevant to our narrative is\ the
\textit{chain rule}. For a time-independent field $F$ that is a restriction of
an ambient field $F$, such as the gravitational potential, the derivative
$\dot{\nabla}F$ is given by the formula%
\begin{equation}
\dot{\nabla}F=C\frac{\partial F}{\partial n}, \label{Chain}%
\end{equation}
where $\partial F/\partial n$ is the normal derivative of $F$. Note that,
given our Tensor Calculus context, $N^{i}\nabla_{i}F$ would have been a more
appropriate way to express this quantity, but we will continue to use
$\partial F/\partial n$ because it is the more familiar symbol to most readers.

With these essential elements in place, we are ready to turn our attention to
the problems at hand. Note that all elements presented here are valid for any
hypersurface, be it a surface in space or a curve in the plane or on a
surface. Furthermore, if we were to omit any mention of the vector normal
$\mathbf{N}$, all presented elements could be extended to higher dimensions,
as well.

\section{Time evolution of integrals}

First consider the volume integral%
\begin{equation}
\int_{\Omega}Fd\Omega
\end{equation}
of an invariant field $F$ defined over the time-dependent domain $\Omega$ with
boundary $S$. Note that we have dropped the explicit dependence of $\Omega$
and $S$ on $t$ since all domains and surfaces are now understood to be
time-dependent. Then the rate of change of this integral is given by the
formula%
\begin{equation}
\frac{d}{dt}\int_{\Omega}Fd\Omega=\int_{\Omega}\frac{\partial F}{\partial
t}d\Omega+\int_{S}CFdS, \label{dV}%
\end{equation}
where $C$ is referred to the \textit{exterior} normal. This formula is easily
seen to be a direct generalization of the Fundamental Theorem of Calculus%
\begin{equation}
\frac{d}{dt}\int_{a}^{b\left(  t\right)  }F\left(  t,x\right)  dx=\int
_{a}^{b\left(  t\right)  }\frac{\partial F\left(  t,x\right)  }{\partial
t}dx+b^{\prime}\left(  t\right)  F\left(  t,b\right)  .
\end{equation}

Next, let us turn our attention to the surface integral%
\begin{equation}
\int_{S}FdS,
\end{equation}
of an invariant field $F$ defined on the hypersurface $S$. The evolution of
this integral is governed by the formula%
\begin{equation}
\frac{d}{dt}\int_{S}FdS=\int_{S}\dot{\nabla}Fds-\int_{S}CFB_{\alpha}^{\alpha
}dS. \label{dS}%
\end{equation}
Note that, unlike the volume integral, the choice of normal here is immaterial
since the sign of the mean curvature $B_{\alpha}^{\alpha}$ also depends on the
choice of normal and therefore the combination $CB_{\alpha}^{\alpha}$ is
independent of it.

We are now ready to show the application of the Calculus of Moving Surfaces to
specific variational problems.

\section{The isoperimetric problem}

Let us begin with the isoperimetric problem. To accommodate the fixed-volume
constraint, form the Lagrangian%
\begin{equation}
L\left(  t\right)  =\int_{S}FdS+\lambda\left(  \int_{\Omega}d\Omega
-V_{0}\right)  ,
\end{equation}
where the time-dependence comes from the virtual evolution $\Omega\left(
t\right)  $ and $S\left(  t\right)  $ of the domain and its boundary.

Applying equations (\ref{dV}) and (\ref{dS}), we find:%
\begin{equation}
L^{\prime}\left(  t\right)  =\int_{S}C\left(  -B_{\alpha}^{\alpha}%
+\lambda\right)  dS.
\end{equation}
Equating $L^{\prime}\left(  0\right)  $ to $0$, we find that the minimal
surface satisfies the equation%
\begin{equation}
\int_{S}C\left(  -B_{\alpha}^{\alpha}+\lambda\right)  dS=0
\end{equation}
for any possible evolution $S\left(  t\right)  $ and thus for any velocity
$C$. Therefore, $C$ may be treated as an independent variation, and thus the
quantity it multiplies in the integrand must vanish, i.e.%
\begin{equation}
B_{\alpha}^{\alpha}=\lambda.
\end{equation}
In other words, the solution is a hypersurface of constant \textit{mean
curvature}, i.e. a circle in two dimensions and a sphere in three dimensions.

If instead we consider a surface patch $S$ not subject to a volume constraint,
but with a fixed boundary, the analysis simplifies further and takes us
straight from the integral%
\begin{equation}
E\left(  t\right)  =\int_{S}dS
\end{equation}
to the equilibrium condition%
\begin{equation}
B_{\alpha}^{\alpha}=0.\tag{%
\ref{BAa = 0}%
}%
\end{equation}
In other words, the \textit{mean curvature vanishes}. When interpreted in
specific coordinate systems, this criterion delivers both Euler's equation
(\ref{Euler}) as well as Lagrange's equation (\ref{Lagrange}). However, the
equilibrium equation (\ref{BAa = 0}), in the form it emerged from the Calculus
of Moving Surfaces, offered a ready geometric interpretation. This is the
simplest illustration of the fact that, in this regard, the Calculus of Moving
Surfaces is superior to the conventional approach.

We now turn our attention to the brachistochrone.

\section{The brachistochrone}

Consider the motion of a material particle, initially at rest, constrained to
slide frictionlessly along a planar curve segment $S$ under the influence of a
general gravitational potential $P$. Assuming, with no loss of generality,
that $P=0$ at the initial location, the conservation of energy reads%
\begin{equation}
P+\frac{1}{2}v^{2}=0.
\end{equation}
Thus, the speed $v$ of the material particle is $\sqrt{-2P}$, and therefore
the total time $T$ it takes to travel along the curve segment is given by%
\begin{equation}
T=\int_{S}\frac{dS}{\sqrt{-2P}}.
\end{equation}

Since the actual concept of time is present in this problem, let us use $\tau$
as the parameter for the virtual motion of the curve $S$. According to
equation (\ref{dS}), the first variation given by the equation%
\begin{equation}
T^{\prime}\left(  \tau\right)  =\int_{S}\left(  \dot{\nabla}\frac{1}%
{\sqrt{-2P}}-\frac{CB_{\alpha}^{\alpha}}{\sqrt{-2P}}\right)  dS.
\end{equation}
According to the chain rule (\ref{Chain}),%
\begin{equation}
\dot{\nabla}\left(  -2P\right)  ^{-1/2}=C\frac{\partial P/\partial n}{\left(
-2P\right)  ^{3/2}}%
\end{equation}
and therefore we find that
\begin{equation}
T^{\prime}\left(  \tau\right)  =\int_{S}C\left(  \frac{\partial P/\partial
n}{\left(  -2P\right)  ^{3/2}}-\frac{B_{\alpha}^{\alpha}}{\left(  -2P\right)
^{1/2}}\right)  dS.
\end{equation}
Once again, since $C$ can be treated as an independent variation, we arrive at
the brachistochrone equation%
\begin{equation}
\frac{\partial P}{\partial n}=-2PB_{\alpha}^{\alpha}.
\end{equation}
For consistency with minimal surfaces, rewrite this equation in the form that
isolates the mean curvature, i.e.%
\begin{equation}
B_{\alpha}^{\alpha}=-\frac{1}{2P}\frac{\partial P}{\partial n}.\label{Brach}%
\end{equation}
It should be noted that this equation remains valid for the problem in which
the material particle is constrained to move on a curved surface, in which
case $B_{\alpha}^{\alpha}$ would be interpreted as the geodesic curvature and
$\partial P/\partial n$ as the conormal derivative.

As with minimal surfaces previously, we have arrived at an equilibrium
equation expressed entirely in terms of geometric quantities. As unlikely as
it may seem, we have not been able to find this form of the brachistochrone
equation in the literature. Meanwhile, given its analogy to the classical
equation (\ref{BAa = 0}) for minimal surfaces, it deserves greater prominence.

It is noteworthy that equation (\ref{Brach}) could be derived from the
so-called \textit{ray equation} \cite{BornWolf1999} in Optics by
differentiating along the curve. In fact, tracing the path of a ray in a
medium with a variable index of refraction was the very basis Johann
Bernoulli's original approach to solving the brachistochrone problem
\cite{Goldstine1980Calculus}, \cite{SandersonStrogatz2016}.\ Nevertheless,
equation (\ref{Brach}) appears to have not been formulated before. If this
indeed proves to be true, it will further illustrate that geometric
interpretation is often outside the scope of variational analyses.

Lastly, we will demonstrate how to reduce the invariant criterion
(\ref{Brach}) to an equation with ordinary derivatives that would have
resulted from a conventional application of the Euler-Lagrange equation. This
step is needed when one wishes to find a specific solution to the problem by
analytical or numerical means.

Suppose that the desired curve is given in standard Cartesian coordinates
$x,y$ by the equations
\begin{align}
x &  =x\left(  \gamma\right)  \\
y &  =y\left(  \gamma\right)  ,
\end{align}
where $\gamma$ denotes the curve coordinate $S^{1}$. Then the gravitational
potential $P$ is given by%
\begin{equation}
P\left(  x,y\right)  =mgy
\end{equation}
where $m$ is the mass of the material particle and $g$ is the acceleration of
gravity. Since the components of the gradient of $P$ are $\left(  0,mg\right)
$, and the components of the normal are $\left(  y_{\gamma},-x_{\gamma
}\right)  /\sqrt{x_{\gamma}^{2}+y_{\gamma}^{2}}$, we have%
\begin{equation}
\partial P/\partial n=-\frac{mgx_{\gamma}}{\sqrt{x_{\gamma}^{2}+y_{\gamma}%
^{2}}}.
\end{equation}
The equation for the mean curvature corresponding to the above choice of
normal is
$\backslash$%
\begin{equation}
B_{\alpha}^{\alpha}=\frac{x_{\gamma}y_{\gamma\gamma}-y_{\gamma}x_{\gamma
\gamma}}{\left(  x_{\gamma}^{2}+y_{\gamma}^{2}\right)  ^{3/2}}.
\end{equation}
Therefore, the equilibrium equation (\ref{Brach}) reads%
\begin{equation}
\frac{x_{\gamma}y_{\gamma\gamma}-y_{\gamma}x_{\gamma\gamma}}{\left(
x_{\gamma}^{2}+y_{\gamma}^{2}\right)  ^{3/2}}=-\frac{x_{\gamma}}%
{2y\sqrt{x_{\gamma}^{2}+y_{\gamma}^{2}}}%
\end{equation}
or, equivalently,%
\begin{equation}
2y\left(  x_{\gamma}y_{\gamma\gamma}-y_{\gamma}x_{\gamma\gamma}\right)
+x_{\gamma}\left(  x_{\gamma}^{2}+y_{\gamma}^{2}\right)  =0
\end{equation}
Finally, when we parameterize the curve by $x$, i.e. $\gamma=x$, we arrive at
the well-known brachistochrone equation%
\begin{equation}
2yy_{xx}+y_{x}^{2}+1=0.
\end{equation}

\section{Conclusion}

Since its introduction, the Euler-Lagrange equation has been at the center of
the Calculus of Variations. It has been deeply studied and the accumulated
body of knowledge in this area is not to be dismissed. However, we have
demonstrated that, by analyzing the problem directly in its geometric form,
the Calculus of Moving Surfaces preserves the connection to the underlying
geometry in a way that Lagrange's analytical approach cannot. Furthermore, in
the tradition of Tensor Calculus, the Calculus of Moving Surfaces keeps the
analytical complexity at bay and otherwise prevents the undesirable artifacts
of coordinate analysis from creeping into the calculation. It therefore
demonstrates \textit{a higher degree of perfection} compared to the
conventional approach and it would therefore serve to greatly benefit the
scientific community if it were more widely adopted.

\bibliographystyle{abbrv}
\bibliography{Classics,Misc,PGrinfeld}

\end{document}